\documentclass{birkjour}

\usepackage{latexsym} 
  \usepackage[all]{xy}
  \usepackage{amsfonts} 
  \usepackage{amsthm} 
  \usepackage{amsmath} 
  \usepackage{amssymb}
  \usepackage{pifont}  
  \usepackage{enumerate}
 \newtheorem{thm}{Theorem}[section]

 \theoremstyle{definition}
 
 \theoremstyle{remark}

 \numberwithin{equation}{section}
 
 \def\bs{\begin{statement}}
\def\es{\end{statement}}
  \swapnumbers
  \newtheorem{statement}[thm]{}

  \newcounter{zlist}

  \newcounter{blist}

  \newcounter{rlist}

\def\ot{\otimes}

\def\AA{{\mathbb A}}

\def\CC{{\mathbb C}}

\def\KK{{\mathbb K}}

\def\RR{{\mathbb R}}

\def\ZZ{{\mathbb Z}}

\def\ii{{\mathfrak I}}

\def\tt{{\mathfrak T}}

\newcommand{\Cc}{\mathcal{C}}

\def\O{{\bf O}}

\def\*C{{}^*\hspace*{-1pt}{\Cc}}

\def\text#1{{\rm {\rm #1}}}

 \def\1{\mathbf{1}}

\begin{document}

%-------------------------------------------------------------------------
% editorial commands: to be inserted by the editorial office
%
%\firstpage{1} \volume{228} \Copyrightyear{2004} \DOI{003-0001}
%\seriesextraline{This is the Concrete Title of this Book\br H.E. R and S.T.C. W, Eds.}
%\firstpage{1}

\title[Differential and integral forms on non-commutative algebras]
 {Differential and integral forms on non-com\-mutative algebras}

%----------Author 1
\author{Tomasz Brzezi\'nski}

\address{%
Department of Mathematics, Swansea University, 
  Swansea SA2 8PP, U.K.\ \newline 
Department of Mathematics, University of Bia{\l}ystok, K.\ Cio{\l}kowskiego  1M,
15-245 Bia\-{\l}ys\-tok, Poland}

\email{T.Brzezinski@swansea.ac.uk}

%\thanks{This work was completed with the support of our
%%----------Author 2
%\author{A Second Author}
%\address{The address of\br
%the second author\br
%sitting somewhere\br
%in the world}
%\email{dont@know.who.knows}
%%----------classification, keywords, date


\date{October 2016}

%\begin{abstract}
%The aim of this work is to provide the contributors to journals or to
%multi-authored books with an easy-to-use and flexible class file compatible
%with \LaTeX\ and \AmS-\LaTeX.
%\end{abstract}

%%% ----------------------------------------------------------------------
\maketitle
%%% ----------------------------------------------------------------------
%\tableofcontents

\section{Introduction}
These lectures describe an algebraic approach to differentiation and integration that is characteristic for non-commutative geometry. The material contained in Section~\ref{sec.diff} is standard and can be found in any text on non-commutative geometry, for example \cite{Dub:lec}.  \ref{item.multi} and \ref{item.free}, which describe concepts introduced in \cite{BrzElK:int}, are exceptions. The bulk of Section~\ref{sec.int} is  based on \cite{Brz:con} and \cite{BrzElK:int}, while the Berezin integral example is taken from \cite{Brz:div}.
\section{Differential forms}\label{sec.diff}
\bs{\bf Differential graded algebras.} 
A {\em differential graded algebra} is a pair $(\Omega, d)$, where $\Omega = \oplus_{n\in \ZZ} \Omega^n$ is a graded algebra and $d: \Omega \to \Omega$ is a degree-one map that squares to zero and satisfies the graded Leibniz rule. Note that $\Omega^0$ is an associative algebra and  all the $\Omega^n$ are $\Omega^0$-bimodules.
\es

\bs{\bf Differential calculus.} 
Given an associative algebra $A$ (over a field $\KK$ of characteristic not 2), by a {\em differential calculus over $A$} we mean  a differential graded algebra $(\Omega A, d)$, such that $\Omega^0 A= A$, $\Omega A$ is generated by $\Omega^1A =Ad(A)$, and $\Omega^n A=0$, for all $n<0$. 
%We say $(\Omega A, d)$ is connected provided $\ker d|_A = \KK$. 
A calculus is said to be {\em $N$-dimensional}, if $\Omega^N A \neq 0$ and $\Omega^n A =0$, for all $n>N$. The pair $(\Omega^1 A, d:A \to \Omega^1 A)$ is called a {\em first-order differential calculus}.
\es

\bs{\bf Universal differential calculus.} Every algebra admits the {\em universal differential calculus}, defined as the tensor product algebra over the kernel of the multiplication map $\mu$ on $A$, with the exterior derivation $d: a\mapsto 1\ot a - a\ot 1$, for all $a\in A$, and then extended to the whole of $T_A(\ker \mu)$ by the graded Leibniz rule. 

Every first-order differential calculus $(\Omega ^1A,d)$ can by extended universally to the full calculus by defining $\Omega A$ as the quotient of the tensor algebra $T_A(\Omega^1A)$ by the relations coming from the graded Leibniz rule and $d^2=0$.
\es

\bs{\bf Volume form.} An $N$-dimensional calculus is said to {\em admit a volume form} if $\Omega^NA \cong A$ as a left and right $A$-module. Any free generator $v$ of $\Omega^NA$ as a left and right $A$-module (if it exists) is called a {\em volume form}. 
%A volume form induces an algebra automorphism $\varphi :A \to A$, by $v\varphi(a) = a v$, and a module isomorphism $\pi: \Omega^NA \to A$, by $v \pi(\omega) = \omega$.
\es

\bs{\bf Skew multi-derivations.}\label{item.multi} If $(\Omega A,d)$ is such that $\Omega^1A$ is a finitely generated as a left $A$-module, then any left $A$-module basis $\{\omega_1, \ldots , \omega_n\}$ of $\Omega^1A$ induces maps $\partial_i, \sigma_{ij}: A\to A$, $i,j=1,\ldots , n$, by
\begin{equation}\label{multi}
d a = \sum_i \partial_i(a) \omega_i, \qquad \omega_i a = \sum_j \sigma_{ij}(a)\omega_j.
\end{equation}
These necessarily satisfy
\begin{subequations}\label{sig.part}
\begin{equation}\label{sig}
\sigma_{ij}(1) =1\qquad   \sigma_{ij}(ab) = \sum_{k} \sigma_{ik}(a)\sigma_{kj}(b),
\end{equation}
\begin{equation}\label{part} 
\partial_j(ab) = \sum_i\partial_i(a)\sigma_{ij}(b) + a\partial_j(b).
\end{equation}
\end{subequations}
A system $(\partial_j, \sigma_{ij})_{i,j=1}^n$ is called a {\em skew multi-derivation}. Any skew  multi-derivation induces a calculus on $A$ by formulae \eqref{multi}, provided there exist $a_i^{\alpha_i}, b_i^{\alpha_i} \in A$ such that $\sum_{\alpha_i} a_i^{\alpha_i} \partial_j(b_i^{\alpha_i} ) = \delta_{ij}$. If such elements exist $(\partial_j, \sigma_{ij})_{i,j=1}^n$ is said to be {\em orthogonal}.

The conditions \eqref{sig} are equivalent to the statement that the map $\mathbf{\sigma}: A\to M_n(A)$, $a\mapsto (\sigma_{ij}(a))_{i,j=1}^n$, where $M_n(A)$ denotes the ring of $n\times n$ matrices with entries from $A$, is an algebra homomorphism.
\es

\bs{\bf Free multi-derivations.} \label{item.free} A skew  multi-derivation $(\partial_j, \sigma_{ij})_{i,j=1}^n$ is said to be {\em free}, provided there exist $\bar{\sigma}_{ij}, \hat\sigma_{ij} :A\to A$ that satisfy conditions \eqref{sig} and are such that, for all $a\in A$,
$$
\sum_k\bar\sigma_{jk}(\sigma_{ik}(a)) \!=\! \sum_k \sigma_{kj}(\bar\sigma_{ki}(a)) \!=\!  \sum_k\hat\sigma_{jk}(\bar\sigma_{ik}(a)) \!=\!  \sum_k \bar\sigma_{kj}(\hat\sigma_{ki}(a)) \!=\! \delta_{ij} a.
$$
If the matrix $\mathbf{\sigma} = (\sigma_{ij})_{i,j=1}^n$ is triangular with invertible diagonal entries, then $(\partial_j, \sigma_{ij})_{i,j=1}^n$  is free.
\es

\bs{\bf Diagonal and skew $q$-derivations.}  If $\mathbf{\sigma} = (\sigma_{ij})_{i,j=1}^n$ is diagonal, then relations \eqref{part} separate into twisted Leibniz rules,
$
\partial_i(ab) = \partial_i(a)\sigma_{ii}(b) + a\partial_i(b),$  for all $i=1, \ldots, n$. Furthermore, if  $\sigma_{ii}$ is an invertible map and there exists $q_i\in \KK$ such that $\sigma_{ii}^{-1}\circ\partial_i\circ \sigma_{ii} = q_i \partial_i$, then $(\partial_i,\sigma_{ii})$ is called a {\em skew $q_i$-derivation}.
\es

%\bs{\bf Skew derivations and calculus.}
%\es

\bs{\bf Calculus for $q$-polynomials.} A {\em $q$-polynomial algebra} or the {\em quantum plane} is the algebra $A = \KK_q[x,y]$ generated by $x,y$ subject to the relation $xy = qyx$. 
%If, in addition, we assume that $x,y$ are invertible, then $A = \KK_q[x^{\pm 1},y^{\pm 1}]$ is called a {\em Laurent $q$-polynomial algebra} (or the {\em quantum trous}). 
The elements of $A$ are finite combinations of monomials $x^ry^s$.

The algebra $A$ admits a first-order calculus freely generated by one-forms $dx, dy$ and relations:
$$
dxx = p xdx, \quad dyx = pq^{-1}xdy, \quad dxy = qydx +(p-1)xdy, \quad dyy = p ydy,
$$
where $p$ is a non-zero scalar. It is understood that $d:x\mapsto dx, y\mapsto dy$. The universal extension of this calculus necessarily yields $dx dy = -qp^{-1}dydx$, $(dx)^2 = (dy)^2 =0$, and it is  a two-dimensional calculus with a volume form, e.g.\  $v = dxdy$. 
%The associated algebra automorphism comes out as
%$$
%\varphi: f(x,y) \mapsto f(q^{-1}s^{-2} ...
%$$

Setting $\omega_1 = dx$, $\omega_2 = dy$ one easily finds that the associated skew multi-derivation $(\partial_i, \sigma_{ij})_{i,j=1}^2$ is free with the upper-triangular matrix-valued endomorphism
$$
\sigma(x^ry^s) = 
\begin{pmatrix} 
p^rq^s x^r y^s & p^r(p^s-1) x^{r+1}y^{s-1} \cr 0 &p^{r+s}q^{-r} x^ry^s
\end{pmatrix} .
$$
\es
\bs{\bf Inner calculus.} \label{Jackson} A calculus $(\Omega A,d)$ is said to be {\em inner} if there exists $\theta \in \Omega^1 A$, such that $d(\omega) = \theta \omega - (-1)^n\omega$, for all $\omega \in \Omega^nA$. Note that $d^2(\omega) = 0$ implies that $\theta^2$ is central in $\Omega A$. Also $\theta$ satisfies the Cartan-Maurer equations $d \theta = 2 \theta^2$.

As an example, consider a one-dimensional calculus on the Laurent polynomial ring $A = \KK[x,x^{-1}]$, given by  Jackson's $q$-derivation 
\begin{equation}\label{Jac}
\partial_q (f) = \frac{f(qx) - f(x)}{(q-1)x},
\end{equation}
where taking the limit is understood in case $q=1$. If $q\neq 1$, this calculus is inner with $\theta = \frac{1}{q-1} x^{-1}dx$, otherwise it is not inner.
\es

\section{Integral forms}\label{sec.int}
\bs{\bf Divergence.}
Let $(\Omega A, d)$ be a differential calculus on an algebra $A$. We will denote by $\ii_nA$, the Abelian group of all right $A$-linear maps $\Omega^n A\to A$. For all $n\geq m$, consider maps
$$
\cdot : \ii_nA \ot \Omega^m A \to \ii_{n-m} A, \qquad f\ot \omega \mapsto f\cdot \omega, \quad (f\cdot \omega)( \omega ' ) =  f(\omega\omega').
$$
In particular, $\cdot$ makes $\ii_n A$ into a right $A$-module. 

A {\em divergence} is a linear map $\nabla_0: \ii_1 A\to A$, such that, for all $a\in A$, $\nabla_0(f\cdot a) = \nabla_0(f) a+ f(da)$. A divergence is extended to a family of maps $\nabla_n:  \ii_{n+1} A\to  \ii_{n} A$, by
$\nabla_n(f) (\omega) = \nabla_0(f\cdot \omega) + (-1)^{n+1}f(d\omega)$, for all $\omega\in \Omega^n A$. 

The cokernel map $\Lambda : A\to  A/{\mathrm{Im}} \nabla_0$ is called the {\em integral} associated to $\nabla_0$.
\es

\bs{\bf Integral forms.}
A divergence $\nabla_0$ is said to be {\em flat}, provided $\nabla_0\circ\nabla_1 =0$. It is then the case that, for all $n$, $\nabla_n\circ \nabla_{n+1} =0$, and hence there is a complex,
$%\begin{equation}\label{int}
\xymatrix{ 
\cdots \ar[r]^-{\nabla_2} & \ii_2 \ar[r]^-{\nabla_1} &  \ii_1 \ar[r]^-{\nabla_0} &  A ,}
$%\end{equation}
known as the complex of {\em integral forms}. 
\es

\bs{\bf The inner case.} If  $(\Omega A, d)$ is an inner differential calculus with the exterior derivative given by a graded commutator with $\theta \in \Omega^1 A$, then $\nabla_0: \ii_1 A\to A$, $f\mapsto -f(\theta)$ is a divergence. One easily finds that $\nabla_1(f)(\omega) = f(\theta \omega)$, and so this divergence is flat, provided $\theta^2=0$. 
\es

\bs{\bf Divergences and multi-derivations.}\label{multi.div} Let $(\partial_i, \sigma_{ij}; \bar\sigma_{ij}, \hat\sigma_{ij})_{i,j =1}^n$ be a free skew multi-derivation on $A$, and let $(\Omega^1A, d)$ be the associated first-order calculus with generators $\omega_1, \ldots ,\omega_n$. Let $\xi_i\in\ii_i$ be the dual basis to the $\omega_i$, i.e.\  the $\xi_i$ are given by $\xi_i(\omega_j) = \delta_{ij}$. 
%Define $\partial_i^\sigma := \sum_{j,k} \bar\sigma_{kj}\circ\partial_j\circ \hat\sigma_{ki}$. 
Then
\begin{equation}\label{nabla}
\nabla_0:  \ii_1 A\to A, \qquad f\mapsto \sum_{i,j,k}  \bar\sigma_{kj}\left(\partial_i\left(f\left(\hat\sigma_{ki}\left(\omega_i\right)\right)\right)\right),
\end{equation}
is a unique divergence such that $\nabla_0(\xi_i) =0$, for all $i=1, \ldots ,n$.

In particular, if $\mathbf{\sigma}$ is diagonal and all the $(\partial_i,\sigma_{ii})$ are skew $q_i$-derivations, then 
$$
\nabla_0(f) = \sum_i q_i \partial_i\left(f\left(\omega_i\right)\right).
$$
\es

\bs{\bf Cauchy's integral formula.} Let $A= \KK[x,x^{-1}]$ be the Laurent polynomial ring with the one-dimensional calculus $(\Omega A, d)$ given by  Jackson's $q$-derivative  \eqref{Jac}. In this case $\partial_q$ is twisted by the automorphism $\sigma(f(x)) = f(qx)$, and $(\partial_q, \sigma)$ is a skew $q$-derivation. $\Omega^1A$ is generated by $dx$, and hence the  corresponding divergence \eqref{nabla} is  $\nabla_0(f) = q\partial_q\left(f\left(dx\right)\right)$. For all $f\in \ii_1 A$, define $f_x\in \KK[x,x^{-1}]$ by
$
f_x(x) := f(dx).$ Then
$$
\nabla_0(f) = q\frac{f_x(qx)-f_x(x)}{(q-1)x}.
$$
The image of $\nabla_0$ consists of all of $\KK[x,x^{-1}]$ except the monomials $\alpha x^{-1}$. Therefore, the integral is
$$
\Lambda: \KK[x,x^{-1}] \to \KK, \qquad a \mapsto \mathrm{res}(a)\Lambda(x^{-1}).
$$
In case $\KK = \CC$ we can normalise $\Lambda$ as $\Lambda(x^{-1}) =  2\pi i$, and obtain the Cauchy integral formula.
\es

\bs{\bf Calculus on quantum groups.} If $A$ is a coordinate algebra of a compact quantum group (over $\CC$), then every left-covariant differential calculus gives rise to a free multi-derivation \cite{Wor:dif}, and hence there is a canonical divergence $\nabla_0$ \eqref{nabla} and the corresponding integral $\Lambda$. Any right integral $\lambda$ on $A$ (the Haar measure) factors uniquely through $\Lambda$, i.e.\ there exists unique $\varphi:  A/{\mathrm{Im}} \nabla_0 \to \CC$ such that $\lambda = \varphi\circ\Lambda$.
\es

\bs{\bf Berezin's integral.} Let $A$ be a superalgebra of  (integrable) real functions on the supercircle $S^{1\mid 1}$. That is, $A$ consists of $a(x,\vartheta) = a^0(x) + a^1(x)\vartheta $,
where $a^i: [0,1]\to \RR$ are (integrable) functions such that $a^i(0) = a^i(1)$ and $\vartheta$ is a Grassmann variable, $\vartheta^2 =0$. The differentiation on $A$ is defined by
$$
\partial_xa(x,\vartheta) := \frac{da^0(x)}{dx} + \frac{da^1(x)}{dx} \vartheta, \qquad \partial_\vartheta a(x,\vartheta)  := a^1(x).
$$
One easily checks that $\partial = (\partial_x, \partial_\vartheta)$ is a free skew multi-derivation with the twisting matrix-valued endomorphism $\mathbf{\sigma}( a(x,\vartheta))= \begin{pmatrix} a(x,\vartheta) & 0 \cr 0 & a(x,-\vartheta)\end{pmatrix}$.
%, where
%$
%\sigma_{xx}(a(x,\vartheta)) = a(x,\vartheta)$, $\sigma_{\vartheta\vartheta}(a(x,\vartheta)) = a(x,-\vartheta).
%$
The calculus $\Omega^1 A$ is freely generated by $dx$ and $d\vartheta$. The corresponding divergence \ref{multi.div} comes out as
$$
\nabla_0(f)(x,\vartheta) = \partial_x f_x(x,\vartheta) - \partial_\vartheta f_\vartheta(x,\vartheta) = \frac{df_x^0(x)}{dx} - f_\vartheta^1(x) + \frac{df_\vartheta^0(x)}{dx}\vartheta,
$$
where  $f_x(x,\vartheta) := f(dx)(x,\vartheta) = f_x^0(x) + f_x^1(x)\vartheta$ and $f_\vartheta(x,\vartheta) := f(d\vartheta)(x,\vartheta) = f_\vartheta^0(x) + f_\vartheta^1(x)\vartheta$. 

If $a(x,\vartheta)$ is purely even, i.e.\ $a(x,\vartheta) = a(x)$, then setting $f_x(x,\vartheta) = 0$ and
$f_\vartheta(x,\vartheta) = -a(x)\vartheta$ we obtain $a(x) = \nabla_0(f)$. Thus, the integral $\Lambda$ vanishes on the even part of $A$. Since $\Lambda$ is the cokernel map of $\nabla_0$, for all $f\in \ii_1 A$,
\[
0 = \Lambda\circ \nabla_0(f) = \Lambda\left(\frac{df_x^0(x)}{dx} - f_\vartheta^1(x) + \frac{df_\vartheta^0(x)}{dx}\vartheta\right) = \Lambda\left(\frac{d}{dx}f_\vartheta^0(x)\vartheta\right).
\]
On the other hand,  $f_\vartheta^0(0) = f_\vartheta^0(1)$, so $\int_0^1 \frac{d}{dx}f_\vartheta^0(x)dx =0$. By the universality of $\Lambda$, $ \Lambda(a^1(x)\vartheta) = \int_0^1 a^1(x)dx$. Therefore,
\[
 \Lambda(a^0(x) + a^1(x)\vartheta) = \int_0^1 a^1(x)dx,
\]
i.e.\ $\Lambda$ is the Berezin integral on the supercircle.

\es

% ------------------------------------------------------------------------

\subsection*{Acknowledgment}
I would like to thank the organizers and participants of the 2016 Bia\l owie\.za School for creating very friendly and stimulating atmosphere.

% ------------------------------------------------------------------------
\end{document}